\documentclass[12pt]{amsart}

\newcommand{\CC}{\mathbb C}
\newcommand{\Z}{\mathbb Z}
\newcommand{\g}{\mathfrak g}
\newcommand{\Conv}{Conv}
\newcommand{\SL}{SL}
\newcommand{\GL}{GL}
\newcommand{\Ker}{\mathop{\rm Ker}\nolimits}

\begin{document}

\title{On integral of exponent of a homogeneous polynomial}
\author{A. V. Stoyanovsky}
\begin{abstract}
We introduce the notion of $G$-hyper\-geo\-metr\-ic function, where $G$ is a complex Lie group. In the case when $G$ is
a complex torus, this notion amounts to the notion of Gelfand's $A$-hyper\-geo\-metr\-ic function. We show that
the integral $\int e^{P(x_1,\ldots,x_n)}dx_1\ldots dx_n$, where $P$ is a homogeneous polynomial, is a
$GL(n)$-hyper\-geo\-metr\-ic function of algebraic $SL(n)$-in\-var\-iants of the polynomial.
\end{abstract}
\thanks{Partially supported by the grant RFBR 10-01-00536.}
\email{alexander.stoyanovsky@gmail.com}
\address{Russian State University of Humanities}

\maketitle

{\bf 1. Introduction.}
The problem of computing the integral of exponent of a homogeneous form of $n$ variables over $n$-dimensional space
has been posed by V.~V.~Dolotin [1]. He also developed the new techniques of invariant theory (diagram techniques,
discriminant theory), which allows one to compute invariants explicitly [2--4]. The paper [5] by I.~M.~Gelfand and
M.~I.~Graev implies that the integral of exponent of a homogeneous form is a generalized hyper\-geo\-metr\-ic function
of the coefficients of the form.
In the paper [6] by A.~Yu.~Morozov and Sh.~R.~Shakirov this integral is calculated explicitly in several particular
cases as a function of algebraic $\SL(n)$-invariants. In these particular cases the integral is a generalized
hyper\-geo\-metr\-ic function, in the sense of I.~M.~Gelfand, of basic invariants. The authors of [6] conjecture that
the integral is always a generalized hyper\-geo\-metr\-ic function of algebraic invariants of a form.
The purpose of the present paper is to discuss this conjecture. Recall that the notion of generalized
hyper\-geo\-metr\-ic function in the sense of Gelfand is related with an action of complex torus.
We generalize this notion to non-Abelian complex Lie groups. We call the obtained functions
by $G$-hyper\-geo\-metr\-ic functions,
where $G$ is a Lie group. In the case when $G$ is a torus, this notion essentially coincides with the notion of
Gelfand's $A$-hyper\-geo\-metr\-ic function. We show that integral of exponent of a form as a function of basic
invariants is always a $\GL(n)$-hyper\-geo\-metr\-ic function. Thus, in order to make the Morozov--Shakirov conjecture
true, one should extend the notion of $A$-hyper\-geo\-metr\-ic function to $G$-hyper\-geo\-metr\-ic function.
There are very interesting open problems of explicit computation of $G$-hyper\-geo\-metr\-ic functions,
their equations and their singularities, in particular, for integral of exponent of a form.
Interesting computations can be found in [12,13]. Other generalizations of hypergeometric functions to Lie groups
are contained in [5,14].

The author thanks V.~V.~Dolotin, A.~Yu.~Morozov, and Sh.~R.~Sha\-kir\-ov for numerous useful discussions.
\medskip

{\bf 2. Preliminaries.}

{\bf Notations.} $k=(k_1,\ldots,k_n)\in\Z^n$, $k_i\ge 0$, $x^k=x_1^{k_1}\ldots x_n^{k_n}$.

$A$ is the set of all $k$ of one and the same homogeneous degree $d=\sum_{i=1}^n k_i$.

$P(a;x)=\sum_{k\in A}a_kx^k$ is a homogeneous polynomial, where $a=(a_k)_{k\in A}\in\CC^A$.

$J_{n|d}(C;a)=\int_C e^{P(a;x)}dx$
is the integral over an $n$-dimensional real contour $C\subset\CC^n$,
on which the function $e^{P(a;x)}$ rapidly decreases at infinity.

$\Delta(A)$ is the discriminant, i.~e. the set of $a\in\CC^A$ for which the projective hypersurface
$P(a;x)=0$ is singular.

For each face $\Gamma$ of the polytop $\Conv(A)$, the convex hull of the set $A$ (including the polytop itself),
denote by $\pi_\Gamma:\CC^A\to\CC^{\Gamma\cap A}$ the projection, $\pi_\Gamma(a)=(a_k)_{k\in\Gamma\cap A}$,
and let $E(\Gamma)=\pi_\Gamma^{-1}(\Delta(\Gamma\cap A))$. Let $E$ be the union of $E(\Gamma)$ for all $\Gamma$.
$E$ is called the {\it principal $A$-dis\-crim\-in\-ant} [9].

\medskip

{\bf Theorem 1} (corollary of [5]).
{\it The integral $J_{n|d}(C;a)$ satisfies the $A$-hyper\-geo\-metr\-ic system of equations} [7].

{\bf Theorem 2}.
{\it The singularities of the function $J_{n|d}(C;a)$ lie on the discriminant $\Delta(A)$.}

{\it Proof}. From Theorem~1, from computation of the characteristic variety of the $A$-hyper\-geo\-metr\-ic system [7,8],
and from the description of the varieties projectively dual to the closures of the torus orbits [9, Introduction],
it follows that the singularities of the function $J_{n|d}(C;a)$ lie on the principal $A$-discriminant $E$.
Since for $\Gamma\ne\Conv(A)$ the variety $E(\Gamma)$ is not $\SL(n)$-invariant and the function $J_{n|d}(C;a)$
is $\SL(n)$-invariant, this implies that actually the singularities lie only on the variety $E(\Conv(A))=\Delta(A)$.
Q.~E.~D.
\medskip

{\bf 3. The $G$-hyper\-geo\-metr\-ic $D$-module.}

In this section we follow the paper [7, \S2], generalizing some of its results to the case of non-Abelian Lie group $G$.

Let $V$ be a finite dimensional representation of the group $G$. Assume that the unit operator on $V$ lies in the
image of the Lie algebra $\g$ of the group $G$ (otherwise let us replace $G$ by $G\times\CC^*$).
Let $W=\overline{G\cdot v}$ be the closure of the $G$-orbit of a nonzero vector $v\in V$.
Assume that this closure contains finitely many $G$-orbits. Let $V^*$ be the dual
$G$-module. Let $\beta:\g\to\CC$ be a character of the Lie algebra.
Consider the $D_V$-module $M=D_V/(D_V\cdot J_W+D_V\cdot(L_{e_i}-\beta(e_i))_{i=1}^{\dim\g}$,
where $J_W$ is the ideal of the variety $W$, $e_i$ is a basis of the Lie algebra $\g$, $L_{e_i}$ is the Lie derivative.
This $D$-module has been considered in the book [10, Theorem 5.2.12] for $W$ smooth, where it is proved that
the characteristic variety $SS(M)$ is contained in the union of conormal bundles $T^*_{W_\alpha}V$ to the
$G$-orbits $W_\alpha$ in $W$.
For $W$ non-smooth the proof is similar ([7]). Hence the $D$-module $M$ is holonomic.

Consider the Fourier transform $FM$ of the $D_V$-module $M$. This $D_{V^*}$-module corresponds to the system of equations
\begin{equation}
\square_j\,\varphi=0,
\end{equation}
\begin{equation}
(L_{e_i}+\beta(e_i)-\chi_0(e_i))\,\varphi=0
\end{equation}
for a function $\varphi(w)$ on the space $V^*$. Here
$\square_j$ is the differential operator with constant coefficients corresponding to a polynomial equation of
the variety $W$ under the Fourier transform; $\chi_0(e_i)$ is the trace of the action of $e_i$ on $V$.
Since the $D_V$-module $M$ is monodromic [11] and $G_1$-equivariant,
where $G_1\subset G$ is the Lie group with the Lie algebra $\g_1=\Ker\beta$, the same is true also for the
$D_{V^*}$-module $FM$. (In particular, $G_1$ contains the commutator $[G,G]$.) Therefore, the characteristic
variety $SS(FM)$ corresponds to $SS(M)$ under the identification $T^*V^*\simeq V\times V^*\simeq T^*V$.
\medskip

{\bf Definition.} {\it The system of equations \emph(1,2\emph) on the space $V^*$
and the corresponding system on the variety of orbits $V^*/G_1$ is called the $G$-hyper\-geo\-metr\-ic system,
and its holomorphic solutions on an open subset in $V^*$ or in $V^*/G_1$ are called $G$-hyper\-geo\-metr\-ic functions.}
\medskip

{\bf Theorem 3.} {\it The $G$-hyper\-geo\-metr\-ic system has a finite number of solutions at the general point.
The singularities of the solutions lie on the union of conic varieties whose projectivization is
projectively dual to the projectivization of the varieties $W_\alpha$,
and has codimension~$1$.}

{\it Proof.} This follows from the holonomicity of the $D$-module $FM$ and from the description of $SS(FM)$ given
above (as in [7]).

{\bf Theorem 4.} {\it The integral $\varphi(a)=J_{n|d}(C;a)$ is a $\GL(n)$-hyper\-geo\-metr\-ic function on the
space $V^*=\CC^A$.}

{\it Proof.} Equation (1) follows from the fact that each polynomial equation of the variety
$W=\overline{G\cdot v}\subset V$, where $V=(\CC^A)^*$, $G=\GL(n)$, $v=(1,1,\ldots,1)$, vanishes also on the variety
$\overline{(\CC^*)^n\cdot v}$, which coincides with $W$ in this case. Hence system of equations (1) is equivalent
to the first part of the $A$-hyper\-geo\-metr\-ic system.
Equation (2) follows from integration by parts. Q.~E.~D.

\end{document}